%% file: master.tex
\input macros.tex
\pageno = 1
\secnum = 0
\input abstract.tex

\input enlarge.tex
\input refs.tex
\bye

%% file: macros.tex
\hyphenation{para-met-rize}

\def\Q{{\bf Q}}

\def\ref#1#2#3#4#5#6{ #1: #2. #3 {\bf #4}, #5 (#6)}

\def\card{\mathop{\rm card}}

\def\Pic{\mathop{\rm Pic}}
\def\mgbar{{\overline{\cal M}_g}}
\def\mgb#1{\overline{\cal M}_{#1}}
\def\mgn{{\cal M}_{g,n}}
\def\mbar#1#2{\overline{\cal M}_{#1,#2}}
\def\mgnbar{\mbar{g}{n}}
\long\def\define{\which{Definition}{\nevermind}\ \rm}
\def\which#1#2#3\par{\smallskip\global\advance \thmnum by 1 \edef#2{#1 \the\secnum .\the\thmnum}{\smallcap #2.}\thinspace{\sl #3}\par}
\def\thm#1{\which{Theorem}{#1}}

\def\prop#1{\which{Proposition}{#1}}
\def\proof{\smallskip{\sl Proof.\ }}
\def\floor#1{\lfloor #1 \rfloor}
\def\ref#1#2#3#4#5#6{ #1: #2. #3 {\bf #4}, #5 (#6)}

 4
\font\largebold=cmbx10 scaled \magstep 3
 2
\font\sit=cmti8
\font\smallcap=cmcsc10
\font\eight=cmr8
\font\eightsl=cmsl8
\def\Q{{\bf Q}}

\def\chapter#1{\vfill\eject{\largebold\centerline{#1}\par}\vskip .2in}
\def\section#1{\medskip \penalty -1000\global\advance \secnum by 1 
{\bf\indent{{\ifnum\the\secnum>0 \the\secnum. \fi}#1. }}\global\thmnum=0}

\overfullrule=0pt
\def\mapsto{\to}
\def\mapsby#1{\smash{\mathop{\mapsto}\limits^{#1}}}

\newcount\secnum
\secnum = -100
\newcount\thmnum
\thmnum = 0
\def\qed{ \ {\vbox{\hrule height 0pt depth .3pt\hbox{\vrule width .3pt height 6pt
\kern 3pt \vrule width .3pt}\hrule height 0pt depth .3pt}}\smallskip}

\headline = {\ifodd\pageno \rightheadline \else \leftheadline \fi}
\footline = {\ifnum\pageno=1 {\hfill 1}\fi}
\hsize = 5.7 in\hoffset = 0.4 in
\vsize = 7.9 in\voffset = 0.5 in
\def\leftheadline{\the\pageno\smallcap \hfil adam logan\hfil}
\def\rightheadline{\smallcap \hfil \ifnum\pageno > 1 kodaira dimension
of moduli spaces \hfil \the\pageno\fi}

%% file: abstract.tex
\vskip 1in
\centerline{\bf THE KODAIRA DIMENSION OF MODULI SPACES}
\centerline{\bf OF CURVES WITH MARKED POINTS, II}

\medskip

\centerline{By \smallcap Adam Logan}

\medskip

\centerline{\hbox to 1.5 in{\hrulefill}}

\medskip

{\narrower {\noindent \eightsl Abstract.}
\eight We sharpen our previous results on the $g$ and $n$ such that
$\mgnbar$ is of general type for some $g$ with $g+1$ prime.}
\bigskip

%% file: enlarge.tex
\section{Introduction}
This short paper gives some small improvements to the results of [L2].
I am grateful to Sean Keel for pointing out a mistake in [L1] which led me
to believe incorrectly that these calculations did not actually produce new
pairs $(g,n)$ with $\mbar{g}{n}$ of general type.

The main idea is to use the map $s: \mbar{g}{n+2} \mapsto \mbar{g+1}{n}$ given
by identifying the last two marked points.  In contrast to other maps and
correspondences between moduli spaces, this one has the advantage
that the image is not entirely contained in interesting divisors.  (It can
hardly be, because the image is precisely $\delta_0$ on $\mgnbar$.)  Therefore,
there are no worries about pulling back effective divisors by this map to get
an effective divisor.  It is possible, but in my view unlikely, that pulling
back other divisors by this map would produce further improvements.  The reason
that this pullback is beneficial is that there is no effective Brill-Noether
divisor on $\mgbar$ when $g+1$ is prime, and pulling back other divisors by
this map would only give less efficient replacements for divisors known to be
present.

\section{Test Curves}
One common way to determine the class of a divisor $D$ on a variety $V$ is
to pull it back by maps $f_i: W_i \mapsto V$.  If the class of sufficiently
many $f_i^*(D)$ can be determined, the class of $D$ on $V$ emerges.
In previous work, we considered maps to $\mgnbar$ for which $f_i^{-1}(D)$ was
contained in a proper subvariety which could be identified specifically.
Here, instead, we take the $W_i$ to be curves considered as subvarieties
of $\mgnbar$ and determine the degree of
the divisors $f_i^*(D)$.  The curves $W_i$ are called test curves.

A priori, this means that we will determine the
class of $D$ only up to numerical equivalence, but in fact numerical and
algebraic equivalence on $\mgnbar$ are the same.  Harer proved this in the
case $n = 0$ and subsequently Arbarello and Cornalba deduced the general
result from this [AC].  The reader may wish to consult
[HM] for some examples of test-curve calculations on $\mgbar$.
Here we recall the Picard group of $\mgnbar$ (see [AC] for proofs):

\define On $\mgnbar$, let $\lambda$
and $\delta_0$ be the pullbacks from $\mgbar$ of the classes of the same name.
Fix $0 \le i \le \floor{g/2}, S \subseteq \{1,2,\dots,n\}$, where $\card S > 1$
if $i = 0$ and $1 \in S$ if $i = g/2$.  We let $\delta_{i;S}$ be the divisor
class on the moduli stack whose value on a family of curves with smooth
general fiber is
the number of curves in the family that have a disconnecting node whose removal
creates a connected component of genus $i$ with precisely the marked points
corresponding to elements of $S$ on it.  Let $\psi_i$ be the divisor class on
the moduli stack which takes the value $-\pi_* (\sigma_i^2)$ on the family 
${\cal X} \, \mapsby {\pi} \, B$
with the $\sigma_i$ as sections.  On $\mbar{g}{1}$, let $\omega$ be the
relative dualizing sheaf for the map $\mbar{g}{1} \mapsto \mgbar$; on
$\mgnbar$, let $\omega_i = \pi^* \omega$, where $\pi$ is the map
$\mgnbar \mapsto \mbar{g}{1}$ that forgets all but the $i$\/th marked point.
(The reader is cautioned that $\omega_i$ is not the relative dualizing
sheaf for the map $\mgnbar \mapsto \mbar{g}{n-1}$ that forgets the $i$\/th
marked point.)

\thm{\picgen} $\Pic_{\hbox{\sit fun}} \mgnbar$, the Picard group of the moduli
stack $\mgnbar$, is free on the following
generators: $\lambda, \delta_0, \psi_i\, (1 \le i \le n),$ and
$\delta_{i;S}\, (0 \le i \le \floor{g/2}), S \subseteq \{1,2,\dots,n\})$.\qed

\define $P$ on $\mbar{g}{2}$ is the divisor $s^*(\delta_0)$, where
$s: \mbar{g}{2} \mapsto \mgb{g+1}$ is the map which takes $(C, p_1, p_2)$
to $C$ with $p_1$ and $p_2$ identified.

\prop{\pull0}
$P = \delta_0 - \omega_1 - \omega_2 -2 \delta_{0;\{1,2\}}
+ \sum_i \delta_{i;\{1\}}$.

\proof For any curve $C$ in
$\mbar{g}{2}$, we have $P \cdot C = \delta_0 \cdot s_* C$.  It is therefore
sufficient (and necessary) to compute the intersections of $s_* C$ with
$\delta_0$, where $C$ runs over a basis of $(\Pic \mbar{g}{2} \otimes \Q)
^\vee$, which we proceed to do.  Of course $P$ is invariant under the action
of $S_2$ on $\mbar{g}{2}$, which saves us some work.

Our first test curve will be a fixed general curve of genus $g-1$ with two
fixed general marked points, attached at a third general point to a pencil of
cubics.  This has intersection number $12$ with $\delta_0$ and $-1$ with
$\delta_1$ (cf.\ [AC]), the other intersection numbers being 0.
Similarly, if we attach a general curve of genus
$g-2$ to a pencil of curves of genus $2$ (obtained by fixing a sextic curve
and a fixed point in the plane and taking the $2$-to-$1$ covers of the
line branched at the six points of intersection of the line with the sextic),
we obtain a curve whose intersection number with $\delta_0$ is $30$, with
$\delta_2$ is $-1$, and with all other standard generators is $0$.
By the push-pull formula, their intersections with $P$ are likewise $12$ and
$30$ respectively.  This suggests---correctly, as it turns out---that the
coefficient of $\delta_0$ in $P$ is $1$.

To prove this, we continue by computing intersections with the other test
curves.  First, fix a curve $C$ of genus $g$ and a point $p_2$ on it, and
let $p_1$ be a variable point.  The intersection of this curve with $P$
is $-2g$, that is, the
sum of the self-intersections of a fiber and the diagonal of $C \times C$,
less $1$ for each intersection point of each of these sections with another.
There being one intersection point, we subtract $2$ (as it affects both
sections).  Meanwhile, the intersection of this with $\delta_{0;1,2}$ is
$1$, with $\omega_1$ is $2g-2$, and with all the other generators is $0$.

If instead $p_1$ varies and $p_2$ is fixed on a curve of genus $i$, attached
to a curve of genus $g-i$ at a fixed point, we get $-2i$ for the intersection
with $P$; as before we get
$2-2i - 1 - 1$, but this time we subtract $1$ for blowing up the diagonal
section where it meets the point of attachment, but promptly add it back
because an additional nondisconnecting node is created.  This implies that
if $p_1$ moves on the component of genus $g-i$, the intersection number is
$-2(g-i)$; checking this directly is a good exercise.  Meanwhile, this curve
has intersection number $1$ with $\delta_{0;1,2}$, $1$ with $\delta_{i;\{1\}}$
and $\delta_{i;\{2\}}$, and $0$ with the others.

Finally, it is necessary to work out the intersections with two other types
of curves: first, fix a curve of genus $i\ (< g-1)$ with two marked points,
and attach it at a fixed point to a moving point on a curve of genus $g-i$;
second, attach this curve at a moving point to a fixed point on a curve of
genus $g-i$.  I claim that both of those, when pushed forward, meet $\delta_0$
with multiplicity $0$.  This is very easy to see in the first case, and only
slightly harder in the second.  In this second case, both sections would have
self-intersection $0$, but they are both blown up at a point, each one creating
an additional nondisconnecting node.  The intersection number is therefore
$2(1-1) = 0$.  The first class of curves meet only $\delta_{g-i;\emptyset}$
of the standard generators; the second meet also $\delta_{g;\{1\}}$ and
$\delta_{g;\{2\}}$, once each, and $\omega_1$ and $\omega_2$, also
once each.

We have now computed the intersection numbers with a basis for the dual space,
and are in a position to determine the coefficients of $P$.  Solving
the equations we have obtained, we first get 
$$w_1 = -1, d_{0;\{1,2\}} = -2.$$  This
quickly leads to the consequences that $d_{i;\{1\}} = 1$ for
all $1 \le i \le g-1$, and $d_{i;\emptyset} = 0$, which are readily verified
to be consistent with all of the equalities computed above.  (Here,
Latin letters indicate the coefficients of the divisors denoted by similar
Greek letters.)  We conclude by noting that by symmetry, $w_2 = -1$ as well.
\qed

The point, then, is that the pullback of the Brill-Noether class on $\mgb{g+1}$
to $\mbar{g}{2}$ by $s$ is effective.  The proposition above allows us to
compute it, as it is immediate that $s^* \lambda = \lambda$ and that
$s^* \delta_i = \delta_{i,\emptyset} + \delta_{i;\{1,2\}}$.  Recalling that
the class of the Brill-Noether class $BN$ on $\mgb{g+1}$ is
$(g+4)\lambda - (g+2)/6 \delta_0 - \sum_i i(g+1-i) \delta_i$, we calculate
that $$s^* BN = (g+4)\lambda - {(g+2) \over 6}(\delta_0 - \omega_1 - \omega_2
- 2\delta_{0;\{1,2\}}) + \sum_i \delta_{i;\{1\}}) - \sum_i i(g+1-i)
(\delta_{i;\emptyset} + \delta_{i-1;\{1,2\}}).$$

\section{New results}
To find improvements to the main theorem of [L2] is now just a matter
of linear algebra.  We recall a definition and theorem from that paper:

\define Let $a_1, \dots, a_n$ be nonnegative integers adding to $g$.
Define the divisor $d_{g;a_1, \dots, a_n}$ on an open subset of $\mgnbar$
to be the reduced divisor supported on the curves and points $C, p_1, \dots,
p_n$ such that the divisor $\sum_i a_i p_i$ moves in a pencil.  The divisor
$D_{g;a_1, \dots, a_n}$ is its extension to $\mgnbar$.

\thm{\coefofd} (Thm.\ 5.4 of [L2].)
With respect to the standard basis of $\Pic_{\hbox{\sit fun}} \mgnbar$,
the coefficient of $\omega_i$ in $D_{g;a_1,\dots,a_n}$ is $a_i(a_i+1)/2$, and
the coefficient of $\delta_{0;\{i,j\}}$ is $-a_i a_j$.  Also, the coefficient
of $\lambda$ is $-1$, and that of $\delta_0$ is $0$.\qed

For example, we show that $\mbar{16}{11}$ has positive
Kodaira dimension.  The Brill-Noether class on $\mgb{17}$ pulls back to
$$20 \lambda + 3 (\omega_1 + \omega_2 - \delta_0
+ 2 \delta_{0;\{1,2\}} - \sum_i \delta_{i;\{1\}}) -
\sum_i i(g-i) (\delta_{i-1;\{1,2\}} + \delta_{i;\emptyset}$$ on
$\mbar{16}{2}$.  Pulling this effective divisor to $\mbar{16}{11}$ in all
different ways and averaging, we get the divisor $$D_1 = 20 \lambda + 
6/11 \sum \omega - 3 \delta_0 + \dots$$ to be effective on $\mbar{16}{11}$.
On the other hand, we know that the divisor $D_{16;2,2,2,2,2,1,1,1,1,1,1}$ is
an effective divisor on $\mbar{16}{11}$.  Averaging it over permutations in
an $11$-cycle and using \coefofd\ 
we find $$D_2 = - \lambda + 21/11 \sum \omega + \dots$$ to be
effective.  It is readily checked that $K_{\mbar{16}{11}} = 2D_1/3 + D_2/3 \,
+ \,$ boundary components is effective.  (Only the higher boundary components
remain, and these are easy, using results of [L2].)  If we averaged over
another $11$-cycle in ${\cal S}_{11}$, we would get a different expression of
$K$ as an effective divisor, so $\kappa_{\mbar{16}{11}} > 0$.

Likewise, with $g = 18$ and $n = 9$, the argument to show that $\mgn$ is
of nonnegative Kodaira dimension proceeds similarly.  Above we found an
effective divisor $D_1$ of the form $$22\lambda - 10 \delta_0/3 + 3 \omega +
\dots$$ on $\mbar{18}{2}$, and it pulls back to give $22\lambda - 10
\delta_0/3 + 2/3 \sum \omega + \dots$.  Also the divisor $D_2 =
D_{18;2,2,2,2,2,2,2,2,2}$  is of the form $-\lambda + 3 \sum \omega + \dots$.
Thus $3 D_1/5 + D_2 = 13\lambda - 2 \delta_0 + \sum \omega + \dots$, and
checking the higher boundary components, we conclude the desired result.
Likewise, this allows us to prove that $\mbar{18}{10}$ is of general type.

It is notable that when we pull back $BN$ from $\mgb{23}$ to $\mbar{22}{2}$,
then pull back to $\mbar{22}{3}$ and average, we get a divisor starting
$26 \lambda - 4\delta_0 + 2\sum \omega$, that is, one in which the relevant
coefficients are multiples of those for the canonical.  In particular, this
tells us that $\mbar{22}{3}$ has nonnegative Kodaira dimension and that
$\mbar{22}{4}$ is of general type.

%% file: refs.tex
\medskip
\noindent{\smallcap
Department of Mathematics, University of California, Berkeley, CA 94703}
\smallskip
\noindent{\sl E-mail:} {\tt adam@math.berkeley.edu}


\smallskip
\frenchspacing
\secnum=-1
\centerline{\smallcap References}
\smallskip

\centerline{\hbox to 0.5 in{\hrulefill}}

\smallskip

\leftskip=20.0pt\parindent=-20.0pt\parskip=10pt

\item{[AC]} E. Arbarello, M. Cornalba, {\sl The Picard groups of the moduli
spaces of curves.} Topology {\bf 26}, 153--171, 1987.

\item{[HM]} J. Harris, I. Morrison, {\sl Moduli of Curves.} GTM 167.
Springer-Verlag, 1998.

\item{[L1]} A. Logan, {\sl Moduli spaces of curves with marked points.}
Harvard University doctoral dissertation, June 1999.

\item{[L2]} A. Logan, {\sl The Kodaira dimension of moduli spaces of curves
with marked points.}  To appear.